\documentclass[11pt,reqno]{amsart}
\usepackage[T1]{fontenc}
\usepackage[margin=1in]{geometry}
\usepackage{amsmath,amssymb,amsthm,mathtools}
\usepackage{microtype}
\usepackage[hypertexnames=false,colorlinks=true,linkcolor=blue,citecolor=blue,urlcolor=blue]{hyperref}
\hypersetup{pdftitle={Sharp complete modified logarithmic Sobolev inequality on spheres},pdfauthor={Li Gao and Long Zhao},pdfkeywords={complete modified logarithmic Sobolev inequality, heat semigroup, sphere, matrix entropy, Fisher information, matrix concentration, quantum memory}}
\numberwithin{equation}{section}
\theoremstyle{plain}
\newtheorem{theorem}{Theorem}
\newtheorem{lemma}[theorem]{Lemma}
\newtheorem{proposition}[theorem]{Proposition}

\theoremstyle{definition}
\newtheorem{remark}[theorem]{Remark}
\newtheorem{example}[theorem]{Example}
\newcommand{\tr}{\operatorname{tr}}

\newcommand{\Ric}{\operatorname{Ric}}
\newcommand{\Hess}{\operatorname{Hess}}
\newcommand{\sa}{\mathrm{sa}}
\newcommand{\HS}{\mathrm{HS}}
\newcommand{\E}{\mathbb E}
\newcommand{\I}{\mathcal I}
\newcommand{\dd}{\,\mathrm{d}\mu}
\title[Sharp CMLSI on spheres]{Sharp modified logarithmic Sobolev inequality for matrix-valued functions on spheres}
\author{Li Gao}
\address{School of Mathematics and Statistics, Wuhan University, Wuhan 430072, China}
\address{Wuhan Institute of Quantum Technology, Wuhan 430075, China}
\email{gao.li@whu.edu.cn}
\author{Long Zhao}
\address{School of Mathematics and Statistics, Wuhan University, Wuhan 430072, China}
\email{zhaolong@whu.edu.cn}
\date{October 27, 2026}

\begin{document}

\begin{abstract}
We prove that for every smooth positive definite matrix-valued function $g$ on the unit sphere $S^d$,
\[
 \int_{S^d}\tr(g(\log g-\log\E g))\dd
 \le\dfrac1{2d}\int_{S^d}\tr((\Delta g)\log g)\dd.
\]
Here $\mu$ is normalized area measure, $\E g$ is the spherical mean,
and $\Delta$ is the positive spherical Laplacian. This complete modified logarithmic Sobolev inequality extends Beckner's sharp logarithmic Sobolev inequality, in its
modified form, to matrix-valued functions independent of the matrix size. 
We show that both entropy and Fisher information decay exponentially under
the heat flow with optimal exponent $2d$, whose proof
combines the Bochner identity with a $4$th-tensor estimate in the matrix-valued setting.
\end{abstract}

\maketitle

\section{Introduction and main results}
In~\cite{Beckner93}, Beckner proved sharp Sobolev inequalities on the unit sphere $S^d$, which imply the following
sharp logarithmic Sobolev inequality
\begin{equation}\label{eq:scalar}
 \int v^2\log v^2\dd-(\int v^2\dd)\log(\int v^2\dd)
 \le\dfrac2d\int |\nabla v|^2\dd,
\end{equation}
where $\mu$ is the normalized spherical measure. The circle case $d=1$  is due to Weissler~\cite{Weissler80}. With $f=v^2>0$, the diffusion chain rule rewrites~\eqref{eq:scalar} as a modified logarithmic Sobolev inequality
\begin{equation}\label{eq:modified}\int f\log f-(\int f\dd)\log(\int f\dd)\dd
 \le\dfrac{1}{2d}\int \frac{|\nabla f|^2}{f}\dd.\end{equation}
 The right hand side $\displaystyle \int \frac{|\nabla f|^2}{f}\dd=\int \langle \nabla f, \nabla \log f\rangle \dd$ is called Fisher information.
In this work, we prove that Beckner's sharp logarithmic Sobolev inequality~\cite{Beckner93}, in the above
modified form \eqref{eq:modified}, extends to matrix-valued functions uniformly for arbitrary matrix size.
Let $M_n=M_n(\mathbb C)$ be the $n\times n$ complex matrices, $M_n^{++}$ be its positive invertible
cone, and let $\tr$ denote the standard matrix trace. 
\begin{theorem}\label{thm:main}
Let $S^d$ be the unit sphere for $d\ge2$. For every $n\ge1$ and smooth matrix-valued function $g:S^d\to M_n^{++}$,
\begin{equation}\label{eq:main}
\int_{S^d}\tr(g(\log g-\log\E g))\dd\le\dfrac1{2d}\int_{S^d}\tr((\Delta g)\log g)\dd,
\end{equation}
where $\mu$ is the normalized spherical measure and $\displaystyle\E g=\int_{S^d}g\dd$.
\end{theorem}
Here $\Delta=-\operatorname{div}\nabla$ is the spherical Laplacian generating the heat semigroup
$H_t=e^{-t\Delta}$. The left hand side of \eqref{eq:main} is the matrix-valued relative entropy, while is the right-hand side is the  matrix-valued Fisher information, which is the negative derivative of relative entropy along the heat semigroup,
\begin{align*}
 D(g\|\E g)
 &=\int_{S^d}\tr(g\log g)\dd-\tr((\E g)\log\E g)
 =\int_{S^d}\tr(g(\log g-\log\E g))\dd,
 \\
 \I(g)&=\int_{S^d}\tr((\Delta g)(\log g-\log\E g))\dd=\int_{S^d}\tr((\Delta g)\log g)\dd=-\left.\frac{d}{dt}D(H_tg\|\E g) \right|_{t= 0},
\end{align*}
where  $H_t$ acts entrywise on $g$. As in the scalar case, the above modified log-Sobolev inequality \eqref{eq:main} bounds the relative entropy $D$ by its dissipation $\I$, which by Gr\"onwall's lemma, is equivalent to exponential convergence to its mean in relative entropy 
\[ D(H_t g\|\E g)\le e^{-2d t} D( g\|\E g)\ . \]
Note that for a matrix-valued function $g$, the mean $\mathbb{E}g$ need not commute with the
function values $g(x),x\in S^d$. The entropy relative to this matrix mean is the
quantity that enters quantum information applications. 

Zhao~\cite{Zhao26} obtains the sharp
matrix-valued $S^1$ inequality with constant $d=1$ by a matrix-valued
Wirtinger inequality and Fisher-information dissipation.
On $S^d$, the complete Bakry--{\'E}mery criterion~\cite[Theorem~1.1]{LJL25}
(see also~\cite[Theorem~4.3]{BGJ22}) gives the constant $d-1$ in our normalization,
while first harmonics bound the optimal constant from above by $d$ through the spectral gap. We attain the sharp upper bound \(d\) uniformly in the matrix size $n$. Note that the matrix-valued version of $L_2$-log-Sobolev inequality \eqref{eq:scalar} fails at even $n=2$ \cite[Proposition 3.15]{GJLL25}.

In the proof, we show the Fisher-information dissipation
estimate~\eqref{eq:bochnerlower},
\begin{align}\label{eq:fisher}
 -\frac12\frac{\mathrm d}{\mathrm dt}\I(H_tg)
 \ge d\I(H_tg) 
\end{align}
which by the Gr\"onwall's Lemma again implies $\I(H_tg) \le e^{-2d t}\I(g) $. Integrating this from $t=0$ to $\infty$ gives Theorem \ref{thm:main}. The proof of Fisher-information dissipation combines the Bochner identity with fourth-derivative tensor estimate
arose from joint convexity of the BKM quantum Fisher information metric.

The rest of the paper is organized as follows.
Section~\ref{sec:hessian} develops the entropy differential calculus
and the $4$-th order tensor estimate.
In Section~\ref{sec:bochner}, we recall the geometric notions, including Bochner identity, and prove the Fisher-information dissipation
estimate.
Section~\ref{sec:applications} is devoted to applications to matrix concentration.

\begin{samepage}
\medskip
\noindent\textbf{Acknowledgment.}
This work is partially supported by the National Natural Science
Foundation of China (grant No.~12401163) and the Department of Science
and Technology of Hubei Province (project Nos.~2025EHA041 and
2025AFA044). The authors acknowledge the use of AI tools (ChatGPT 6.1-Sol) in the
exploratory work and in drafting and revising the exposition.
The authors take full responsibility for the mathematical content and
the final manuscript.
\end{samepage}

\section{Entropy derivatives and tensor estimates}\label{sec:hessian}
We start with a discussion on the derivatives of entropy functionals and the estimates of the tensor form that arose from the higher order derivatives.
Let $M_n$ be the algebra of $n\times n$ complex matrices and $(M_n)_{\sa}$ for Hermitian matrices.
Fix $n\ge1$ and set $\Phi(A)=\tr(A\log A)$ as the entropy functional on $M_n^{++}$. Recall the directional derivative of logarithm.
\[ D\log_A[X]=\frac{d}{dt}\log (A+tX)\Big|_{t=0}= \int_0^\infty R_s(A)XR_s(A)\,\mathrm{d}s,\qquad R_s(A)=(A+sI)^{-1},\]
where $I$ is the identity matrix. Denote the symmetric forms
\begin{equation*} h_A=D^2\Phi(A),\qquad T_A=D^3\Phi(A),\qquad Q_A=D^4\Phi(A)
\end{equation*}
for the higher order derivatives. Using the derivative of resolvent
$$DR_s(A)[X]=-R_s(A)XR_s(A)$$ 
and denoting $R_s=R_s(A)$, for Hermitian directions $X,Y,Z,W$, 
\begin{align*} h_A(X,Y)=&\frac{\partial^2 }{\partial t_1 \partial t_2}  \Phi(A+t_1X+t_2Y)\Big|_{t_1=t_2=0}=\int_0^\infty
  \tr\bigl(XR_sYR_s\bigr)\,\mathrm{d}s \\
T_A(X,Y,Z)=& \frac{\partial^3 }{\partial t_1 \partial t_2 \partial t_3}  \Phi(A+t_1X+t_2Y+t_3Z)\Big|_{t_1=t_2=t_3=0}
\\ =& -\int_0^\infty
  \tr\bigl(XR_sYR_sZR_s+XR_sZR_sYR_s\bigr)\,\mathrm{d}s.\\
  Q_A(X,Y,Z,W)=& \frac{\partial^4 }{\partial t_1 \partial t_2 \partial t_3 \partial t_4}  \Phi(A+t_1X+t_2Y+t_3Z+t_4W)\Big|_{t_1=t_2=t_3=t_4=0}
\\=&
\int_0^\infty
\sum_{\pi\in S_3}
\tr \bigl(
W R_s X_{\pi(1)} R_s X_{\pi(2)} R_s X_{\pi(3)} R_s
\bigr)\,ds,
\end{align*}
where \((X_1,X_2,X_3)=(X,Y,Z)\) and the summation is over the permutation group $S_3$. Cyclicity of the trace yields symmetry in all above forms.

These symmetric forms satisfy
\begin{align*}
 Dh_A[Z](X,Y)
 &=\left.\frac{\mathrm d}{\mathrm dt}h_{A+tZ}(X,Y)\right|_{t=0}
 =T_A(X,Y,Z),\\
 DT_A[W](X,Y,Z)
 &=\left.\frac{\mathrm d}{\mathrm dt}T_{A+tW}(X,Y,Z)\right|_{t=0}
 =Q_A(X,Y,Z,W).
\end{align*}
Recall that the inverse logarithmic derivative is given by~\cite{Lieb73}
\[
 (D\log)_A^{-1}(Z)=\int_0^1A^tZA^{1-t}\,\mathrm{d}t.
\]
The second derivative $h_A(X,Y)=\tr(X(D\log)_A(Y))$ is the BKM
inner product. Representing the third derivative in this inner product, we have
\begin{equation}\label{eq:c}
 h_A(c_A(X,Y),Z)=T_A(X,Y,Z),
\end{equation}
where the defining relation of $c_A$ is $(D\log)_A(c_A(X,Y))=(D^2\log)_A(X,Y)$, or more explicitly
\[
 c_A(X,Y)=-(D\log)_A^{-1}\left(
  \int_0^\infty R_s(XR_sY+YR_sX)R_s\,\mathrm{d}s\right).
\]
One analytic input is joint convexity of $h_A$ for the base point $A$ and the
direction $X$ (see \cite[Theorem~3]{Lieb73} and also~\cite{CT14,CH16}). 
\begin{proposition}\label{prop:matrixconvex}
The map $(A,X)\mapsto h_A(X,X)$ is jointly convex on
$M_n^{++}\times(M_n)_{\sa}$.
\end{proposition}
In the following, we write $\|X\|_A^2=h_A(X,X)$ for the norm.
Joint convexity yields the following fourth-order estimate, which is the
logarithmic case of the matrix entropy differential
criterion~\cite[Theorem~3.3(e)]{CH16}.
\begin{lemma}\label{lem:differential}
For every $A\in M_n^{++}$ and $X,Y\in(M_n)_{\sa}$,
\begin{equation}\label{eq:fourth}
 Q_A(X,X,Y,Y)\ge2\|c_A(X,Y)\|_A^2\ge0.
\end{equation}
\end{lemma}

\begin{proof}
Fix $c_A(X,Y)$ fix at the base point $A$. Apply the joint convexity
in Proposition~\ref{prop:matrixconvex} to the affine path
$(A+tX,Y-tc_A(X,Y))$. The resulting scalar function has a nonnegative
second derivative at $t=0$. By~\eqref{eq:c},
\begin{align*}
 0&\le Q_A(X,X,Y,Y)-4T_A(X,Y,c_A(X,Y))+2\|c_A(X,Y)\|_A^2\\
  &=Q_A(X,X,Y,Y)-2\|c_A(X,Y)\|_A^2.\qedhere
\end{align*}
\end{proof}

The preceding bound controls each pair $(X,Y)$ of Hermitian directions. The next lemma
combines the bounds using symmetry of the fourth derivative tensor. Removing
the tensor trace gives the factor needed for the Hessian trace in the
Bochner identity.
\begin{lemma}\label{lem:tensor}
Let $A\in M_n^{++}, k\geq1$ and
$X_1,\ldots,X_k\in(M_n)_{\sa}$. Put
\[
 c_{ij}=c_A(X_i,X_j),\qquad c_*=\sum_{i=1}^kc_A(X_i,X_i),\qquad
 c^0_{ij}=c_{ij}-\dfrac{\delta_{ij}}k c_*.
\]
Then
\begin{equation}\label{eq:tensor}
 \sum_{i,j=1}^k\|c^0_{ij}\|_A^2
 \le\dfrac{k-1}{2k}\sum_{i,j=1}^kQ_A(X_i,X_i,X_j,X_j).
\end{equation}
\end{lemma}

\begin{proof}
First estimate the diagonal part. Put
$q_{ij}=Q_a(X_i,X_i,X_j,X_j)$ and apply~\eqref{eq:fourth} to
$X_i+X_j$ and $X_i-X_j$. Symmetry of $c_a$ gives
$c_a(X_i+X_j,X_i-X_j)=c_{ii}-c_{jj}$.  Full symmetry of the fourth derivative $Q$
cancels the odd terms, so
\[
 q_{ii}+q_{jj}-2q_{ij}\ge2\|c_{ii}-c_{jj}\|_A^2.
\]
The diagonal trace-free norm is related to pairwise differences by
\[
 \sum_{i<j}\|c_{ii}-c_{jj}\|_A^2
 =k\sum_i\|c_{ii}\|_A^2-\|c_*\|_A^2
 =k\sum_i\left\|c_{ii}-\dfrac{c_*}k\right\|_A^2,
\]
so summing the preceding bound over $i<j$ yields
\[
 2k\sum_i\|c^0_{ii}\|_A^2
 \le(k-1)\sum_iq_{ii}-2\sum_{i<j}q_{ij}.
\]
For $i\ne j$, $c^0_{ij}=c_{ij}=c_{ji}$ and
$2\|c_{ij}\|_A^2\le q_{ij}$ by Lemma~\ref{lem:differential}.
Adding the off-diagonal contribution to the diagonal bound gives
\begin{align*}
 2k\sum_{i,j}\|c^0_{ij}\|_A^2
 &=2k\sum_i\|c^0_{ii}\|_A^2+4k\sum_{i<j}\|c_{ij}\|_A^2\\
 &\le(k-1)\sum_iq_{ii}+2(k-1)\sum_{i<j}q_{ij}
  =(k-1)\sum_{i,j}q_{ij}.
\end{align*}
Dividing by $2k$ proves~\eqref{eq:tensor}.
\end{proof}

\section{Bochner identity and Fisher dissipation}\label{sec:bochner}

\subsection{Geometric preliminaries}\label{subsec:geometry}

Let $S^d=\{x\in\mathbb R^{d+1}:|x|=1\}$ be the unit $d$-sphere.
At $x\in S^d$, the tangent space is $T_xS^d=x^\perp$, equipped
with the metric induced by the Euclidean inner product
$\langle\cdot,\cdot\rangle$. Let $\mu$ be the normalized sphere
measure. We recall the geometric constructions used below;
see~\cite[Chapters~2,4,5,7]{Lee18}.

A local orthonormal frame $(e_i)_{i=1}^d$ is a smooth choice of
orthonormal bases for the tangent spaces on an open subset of $S^d$.
For a smooth real function $f$, write
\[
 \partial_i f:=e_i(f).
\]
Thus $\partial_i$ denotes differentiation along $e_i$;
these operators need not be coordinate partial derivatives.
More generally, $Xf$ denotes the directional derivative along a
smooth tangent vector field $X$.

A smooth real one-form $\alpha$ assigns a linear functional
$\alpha_x:T_xS^d\to\mathbb R$ to each point $x$, with smooth
local components $\alpha_i=\alpha(e_i)$. The differential of a
smooth function is the one-form $\mathrm df$ given by
$\mathrm df(X)=Xf$. The metric identifies one-forms with vector
fields: the metric dual $\alpha^\sharp$ is characterized by
\[
 \langle\alpha^\sharp,X\rangle=\alpha(X),\qquad
 \langle\operatorname{grad}f,X\rangle=\mathrm df(X)=Xf.
\]
In particular, $\operatorname{grad}f=(\mathrm df)^\sharp$, and
\[
 \alpha^\sharp=\sum_i\alpha_i e_i,\qquad
 \operatorname{grad}f=\sum_i(\partial_i f)e_i.
\]

To differentiate tangent vector fields, we use the Levi--Civita
connection $\nabla$. It differentiates a vector field $Y$ along $X$
intrinsically and is the unique torsion-free connection compatible
with the metric. Its product rules are
\[
\begin{aligned}
 \nabla_{fX}Y&=f\nabla_XY,\qquad
 \nabla_X(fY)=(Xf)Y+f\nabla_XY,\\
 X\langle Y,Z\rangle
 &=\langle\nabla_XY,Z\rangle+\langle Y,\nabla_XZ\rangle.
\end{aligned}
\]
The induced covariant derivative of a one-form is determined by
\[
 (\nabla_X\alpha)(Y)
 :=X\bigl(\alpha(Y)\bigr)-\alpha(\nabla_XY).
\]
For tensor fields, write $\nabla_i:=\nabla_{e_i}$. In components,
\[
 (\nabla_i\alpha)_j
 :=(\nabla_{e_i}\alpha)(e_j)
 =\partial_i\alpha_j-\alpha(\nabla_{e_i}e_j).
\]
The second term accounts for the variation of the frame.
We also write $\nabla_i\alpha_j$ for $(\nabla_i\alpha)_j$.
In this section, $\operatorname{grad}$ denotes the gradient,
$\partial_i$ differentiates functions, and $\nabla_i$ differentiates
tensor fields covariantly; see~\cite[Propositions~4.15 and~5.17]{Lee18}
for the induced tensor connection and its compatibility with metric contraction.

Applying the covariant derivative to $\mathrm df$ gives the Hessian
$\Hess f=\nabla\mathrm df$. Metric compatibility yields
\[
 \Hess f(X,Y)
 =\langle\nabla_X\operatorname{grad}f,Y\rangle
 =X(Yf)-(\nabla_XY)f.
\]
The torsion-free property makes the Hessian symmetric.
In the local frame,
\[
 \Hess f(e_i,e_j)
 =\partial_i\partial_jf-(\nabla_{e_i}e_j)f.
\]
The metric also defines the pointwise norms
\[
 |\operatorname{grad}f|^2=\sum_i(\partial_i f)^2,\qquad
 |\Hess f|^2=\sum_{i,j}\bigl(\Hess f(e_i,e_j)\bigr)^2.
\]

The divergence of a vector field is the trace of its covariant
derivative. We use the positive Laplace--Beltrami operator:
\[
 \operatorname{div}X
 :=\sum_i\langle\nabla_{e_i}X,e_i\rangle,\qquad
 \Delta f:=-\operatorname{div}(\operatorname{grad}f)
          =-\sum_i\Hess f(e_i,e_i).
\]
These summations are independent of the local orthonormal frame.
For a one-form, metric compatibility gives
\[
 \operatorname{div}(\alpha^\sharp)
 =\sum_i\bigl(\partial_i\alpha_i-\alpha(\nabla_{e_i}e_i)\bigr)
 =\sum_i (\nabla_i\alpha)_i.
\]
Since $S^d$ has no boundary, the divergence
theorem~\cite[Problem~2-22(a)]{Lee18} yields
\begin{equation}\label{eq:stokes}
 \int_{S^d}\sum_i\nabla_i\alpha_i\dd=\int_{S^d}\operatorname{div}(\alpha^\sharp)\dd=0.
\end{equation}
In particular, for smooth functions $f,q$,
\[
 \int_{S^d}f\Delta q\dd
 =\int_{S^d}\langle\operatorname{grad}f,
                         \operatorname{grad}q\rangle\dd,
 \qquad \int_{S^d}\Delta q\dd=0.
\]
Taking $q=f$ shows that $\Delta$ is nonnegative on $L_2(\mu)$.
The associated heat semigroup $H_t=e^{-t\Delta}$ is Markov and
satisfies $\partial_tH_tf=-\Delta H_tf$.

Curvature measures the failure of covariant derivatives to commute:
\[
 R(X,Y)Z
 :=\nabla_X\nabla_YZ-\nabla_Y\nabla_XZ-\nabla_{[X,Y]}Z,
\]
where the Lie bracket satisfies $[X,Y]f=X(Yf)-Y(Xf)$.
The torsion-free property means
$\nabla_XY-\nabla_YX=[X,Y]$.
The Ricci tensor is the contraction
\[
 \Ric(X,Y):=\sum_i\langle R(e_i,X)Y,e_i\rangle.
\]
The unit sphere has constant sectional curvature $1$, hence
\begin{equation}\label{eq:ricci}
 \Ric(X,Y)=(d-1)\langle X,Y\rangle,\qquad
 \Ric_{ij}:=\Ric(e_i,e_j)=(d-1)\delta_{ij};
\end{equation}
see~\cite[Theorem~8.34 and Proposition~8.36]{Lee18}.

For a covariant two-tensor $K$, the induced derivative is
\[
 (\nabla_XK)(Y,Z)
 =X\bigl(K(Y,Z)\bigr)
  -K(\nabla_XY,Z)-K(Y,\nabla_XZ).
\]
We use the component notation
$(\nabla_iK)_{jk}:=(\nabla_{e_i}K)(e_j,e_k)$.
Applied to the Hessian of a smooth scalar function $g$, the classical
Ricci identity~\cite[Theorem~7.14]{Lee18} gives
\begin{equation}\label{eq:hessian-commutation}
 (\nabla_X\Hess g)(Y,Z)-(\nabla_Y\Hess g)(X,Z)
 =\langle R(X,Y)\operatorname{grad}g,Z\rangle.
\end{equation}
Contracting~\eqref{eq:hessian-commutation} in the local orthonormal
frame and using the symmetry of the Hessian gives the contracted
Ricci commutation identity:
\begin{align}\label{eq:bochner}
 \sum_i (\nabla_i \Hess g)_{ij}
 &=\sum_i (\nabla_j \Hess g)_{ii}
   +\sum_k\Ric_{jk} \partial_k g \nonumber\\
 &=-\partial_j\Delta g+\sum_k\Ric_{jk}\partial_k g.
\end{align}

Here $(\nabla_i \Hess g)_{jk}$ denotes the covariant derivative of the
Hessian tensor, and covariant differentiation commutes with taking
its trace.

\subsection{Matrix-valued Fisher information dissipation}

Let $d\ge2$ and let $g:S^d\to M_n^{++}$ be a smooth positive
definite matrix-valued function. The gradient, Laplacian $\Delta$,
and heat semigroup $H_t=e^{-t\Delta}$ act entrywise. Covariant
derivatives act on tangent indices as above and on matrix
coefficients entrywise.
We combine the contracted Ricci identity~\eqref{eq:bochner}, the Ricci
curvature bound~\eqref{eq:ricci}, and the matrix trace-free tensor
estimate in Lemma~\ref{lem:tensor} to prove
\[
 -\frac12\frac{\mathrm{d}}{\mathrm{d}t}\I(H_tg)
 \ge d\I(H_tg).
\]
Integrating this estimate along the heat flow yields the entropy
inequality in Theorem~\ref{thm:main}. The classical counterpart is
the finite-dimensional Bakry--{\'E}mery entropy
method~\cite[pp.~201--203]{BE85}.

We distinguish the Riemannian metric on the sphere from the BKM
inner product on the matrix coefficients. Write
$\langle\cdot,\cdot\rangle_{S^d}$ for the spherical metric.
All matrix coefficients below are Hermitian.
In a local orthonormal frame $(e_i)_{i=1}^d$, write matrix-valued
tangent fields as
\[
 X=\sum_i e_i\otimes X_i,\qquad
 Y=\sum_i e_i\otimes Y_i.
\]
The spherical metric and the Hilbert--Schmidt inner product induce
the pointwise pairing
\[
 \langle X,Y\rangle_{\HS}
 =\sum_{i,j}\langle e_i,e_j\rangle_{S^d}\tr(X_iY_j)
 =\sum_i\tr(X_iY_i).
\]

At each point $x\in S^d$, put $R_s(g)=(g+sI)^{-1}$.
The BKM inner product on the matrix coefficients and its associated
norm are
\[
 h_g(A,B)
 =\tr\bigl(A(D\log)_g(B)\bigr)
 =\int_0^\infty
   \tr\bigl(AR_s(g)BR_s(g)\bigr)\,\mathrm{d}s,
 \qquad
 \|A\|_g^2=h_g(A,A).
\]
Here $g$ abbreviates $g(x)$, so $\|A\|_g$ is a pointwise matrix norm.
Combining $h_g$ with the spherical metric gives
\[
\begin{aligned}
 \langle X,Y\rangle_g
 &=\sum_{i,j}\langle e_i,e_j\rangle_{S^d}\,h_g(X_i,Y_j)
  =\sum_i h_g(X_i,Y_i)\\
 &=\int_0^\infty
   \langle X,R_s(g)YR_s(g)\rangle_{\HS}\,\mathrm{d}s,
\end{aligned}
\]
where the resolvents act on the matrix coefficients of $Y$.
Both pairings are independent of the local orthonormal frame.
Thus the spherical metric contracts the tangent indices,
while $h_g$ pairs the matrix coefficients.

The Hessian components satisfy
\[
 (\Hess g)_{ij}=(\Hess g)(e_i,e_j)
       =(\nabla_i\mathrm dg)(e_j)
       =\partial_i\partial_jg-(\nabla_{e_i}e_j)g.
\]
Here $\mathrm dg(e_j)=\partial_jg$, and $\nabla_i\mathrm dg$
is the covariant derivative of the matrix-valued one-form
$\mathrm dg$. Consequently,
\[
 \operatorname{grad}g=\sum_i e_i\otimes\partial_i g,\qquad
 \sum_i (\Hess g)_{ii}=-\Delta g,\qquad
 \|\Hess g\|_g^2=\sum_{i,j}h_g((\Hess g)_{ij},(\Hess g)_{ij}).
\]
All pairings with tangent vectors in these identities use only
the spherical metric, with matrix coefficients treated entrywise.

Recall that $h_g=D^2\Phi(g)$, $T_g=D^3\Phi(g)$, and
$Q_g=D^4\Phi(g)$, where $\Phi(g)=\tr(g\log g)$.
The dependence of the matrix inner product on $g$ is expressed by
the following product rules. For smooth Hermitian matrix-valued
functions $A,B,C$ and a spatial directional derivative or time
derivative $\partial$,
\begin{align}
 \partial\bigl(h_g(A,B)\bigr)
 &=T_g(\partial g,A,B)
   +h_g(\partial A,B)+h_g(A,\partial B), \label{eq:h}\\
 \partial\bigl(T_g(A,B,C)\bigr)
 &=Q_g(\partial g,A,B,C)
   +T_g(\partial A,B,C)+T_g(A,\partial B,C)+T_g(A,B,\partial C).\label{eq:T}
\end{align}

The same rules hold for tensor-valued arguments, with covariant
derivatives acting on all tangent indices.
Along $g_t=H_tg$, we use $\partial_tg_t=-\Delta g_t$.

The chain rule gives
$\partial_i\log g=(D\log)_g(\partial_i g)$.
Hence integration by parts yields
\begin{equation}\label{eq:Igradient}
\begin{aligned}
 \I(g)
 &=\int_{S^d}\tr\bigl((\Delta g)\log g\bigr)\dd=\int_{S^d}
   \langle\operatorname{grad}g,
          \operatorname{grad}(\log g)\rangle_{\HS}\dd\\
 &=\int_{S^d}
   \langle\operatorname{grad}g,\operatorname{grad}g\rangle_g\dd
  =\int_{S^d}\sum_i\|\partial_i g\|_g^2\dd.
\end{aligned}
\end{equation}

With this notation, the weighted Bochner identity takes the
following form.

\begin{proposition}\label{prop:bochner}
For every smooth $g:S^d\to M_n^{++}$,
\begin{equation}\label{eq:weightedbochner}
\begin{aligned}
 \int_{S^d}\|\Delta g\|_g^2\dd
 &={}(d-1)\I(g)\\
 &+\sum_{i,j}\Bigl(\int_{S^d}
   \|(\Hess g)_{ij}\|_g^2
   +3T_g((\Hess g)_{ij},\partial_i g,\partial_j g)+
   Q_g(\partial_i g,\partial_i g,
       \partial_j g,\partial_j g)\dd\Bigr).
\end{aligned}
\end{equation}
where the matrix-valued Fisher information is 
\begin{align*}
 \I(g)=\int_{S^d}\sum_i\|\partial_i g\|_g^2\dd.
\end{align*}
\end{proposition}

\begin{proof}
Define the real one-form $\eta$ by
\[
 \eta(X)
 :=\sum_j T_g\bigl(\mathrm dg(X),\partial_jg,\partial_jg\bigr),
\]
where $\mathrm dg(X)=Xg
=\langle\operatorname{grad}g,X\rangle_{S^d}$.
The summation over the last two tangent indices is independent
of the local orthonormal frame. Its components are
\[
 \eta_i=\eta(e_i)
 =\sum_jT_g(\partial_i g,\partial_j g,\partial_j g).
\]
Recall that 
\[
 (\Hess g)_{ij}=(\nabla_{e_i}\mathrm dg)(e_j)
 =\partial_i\partial_jg-(\nabla_{e_i}e_j)g.
\]
Using $\operatorname{div}(\eta^\sharp)=\sum_i(\nabla_i\eta)_i$
and the product rule~\eqref{eq:T}, we obtain
\[
 \sum_i(\nabla_i\eta)_i
 =\sum_{i,j}Q_g(\partial_i g,\partial_i g,\partial_j g,\partial_j g)
  -\sum_iT_g(\Delta g,\partial_i g,\partial_i g)
  +2\sum_{i,j}T_g((\Hess g)_{ij},\partial_i g,\partial_j g).
\]
Applying the divergence theorem~\eqref{eq:stokes}, the global integral of the left hand side of the above vanishes, which yields 
\begin{equation}\label{eq:Tdiv}
 \int_{S^d}\sum_i
 T_g(\Delta g,\partial_i g,\partial_i g)\dd
 =\int_{S^d}\sum_{i,j}\left(2T_g((\Hess g)_{ij},\partial_i g,\partial_j g)+
 Q_g(\partial_i g,\partial_i g,
     \partial_j g,\partial_j g)\right)\dd.
\end{equation}

For the weighted Bochner identity, consider the real one-form
$$\omega(X)=h_g(\Delta g,dg(X))
  +\sum_jh_g(\partial_j g,\langle \nabla_{X}\operatorname{grad} g,e_j\rangle)$$ 
whose local components are
\[
 \omega_i=\omega(e_i)
 =h_g(\Delta g,\partial_i g)
  +\sum_jh_g(\partial_j g,(\Hess g)_{ij}).
\]
Applying the Ricci commutation formula \eqref{eq:bochner} entrywise  and the product rules gives
\[
\begin{aligned}
 \operatorname{div}(\omega^\sharp)=\sum_i(\nabla_i\omega)_i
 ={}&\|\Hess g\|_g^2-\|\Delta g\|_g^2
      +\sum_{i,j}\Ric_{ij}h_g(\partial_i g,\partial_j g)\\
 &+\sum_{i,j}T_g((\Hess g)_{ij},\partial_i g,\partial_j g)
      +\sum_iT_g(\Delta g,\partial_i g,\partial_i g).
\end{aligned}
\]
Applying the divergence theorem~\eqref{eq:stokes} again for $\omega$, using the constant Ricci curvature
$\Ric_{ij}:=\Ric(e_i,e_j)=(d-1)\delta_{ij}$, and
substituting~\eqref{eq:Tdiv} proves~\eqref{eq:weightedbochner}.
\end{proof}

We now separate the trace and trace-free parts in the tangent indices. Set
\[
 c_{ij}=c_g(\partial_i g,\partial_j g),\qquad
 a^0_{ij}=(\Hess g)_{ij}+\dfrac{\delta_{ij}}d\Delta g,\qquad
 c^0_{ij}=c_{ij}-\dfrac{\delta_{ij}}d\sum_kc_{kk}.
\]
Since $c_g$ is symmetric and bilinear, the formula
\[
 c(X,Y):=c_g(\mathrm dg(X),\mathrm dg(Y))
\]
defines a symmetric matrix-valued covariant $2$-tensor with components
$c_{ij}$. The trace-free tensors $a^0,c^0$ and their complete
contractions are therefore independent of the local orthonormal frame. For the scalar counterpart of the trace-free Hessian argument,
see~\cite[Section~3, proof of Lemma~3.1]{DE20}.

\begin{theorem}\label{cor:fisher-dissipation}
Let $d\ge2$. For every $n\ge1$ and smooth $g:S^d\to M_n^{++}$, the heat flow
$g_t=H_tg$ satisfies
\begin{equation}\label{eq:bochnerlower}
 -\left.\dfrac12\dfrac{\mathrm{d}}{\mathrm{d}t}\I(g_t)\right|_{t=0^+}-d\I(g)
 \ge\sum_{i,j}\left(\int_{S^d}
      \dfrac d{d-1}\|a^0_{ij}+\beta_dc^0_{ij}\|_{g}^2
       +\gamma_d Q_g (\partial_i g,\partial_i g,\partial_j g,\partial_j g)\dd\right)
 \ge0,
\end{equation}
where the constants $\beta_d$ and $\gamma_d$ depend only on $d$,
\[\beta_d=\frac{2d+1}{2(d+2)},\qquad \gamma_d=\frac{4d-1}{8(d+2)^2}\ .\]
\end{theorem}

\begin{proof}
Write $g=g_t$ and choose a local orthonormal frame independent of
$t$. The heat flow preserves a strictly positive lower bound for $g$, so
we may differentiate the expressions below under the integral.
For this proof, set
\begin{equation}\label{eq:CW}
 \mathcal C=\int_{S^d}\sum_{i,j}T_g((\Hess g)_{ij},\partial_i g,\partial_j g)\dd,
 \qquad
 \mathcal W=\int_{S^d}\sum_{i,j}Q_g(\partial_i g,\partial_i g,\partial_j g,\partial_j g)\dd.
\end{equation}
Since $\partial_tg=-\Delta g$, we have
$\partial_t(\partial_i g)=-\partial_i\Delta g$.
Differentiating~\eqref{eq:Igradient} and using~\eqref{eq:h} gives
\[
\begin{aligned}
 \left.\frac{\mathrm d}{\mathrm dt}\I(g_t)\right|_{t=0^+}
 ={}&\int_{S^d}\sum_iT_g(\partial_tg,\partial_i g,\partial_i g)\dd
     +2\int_{S^d}\sum_i h_g(\partial_t(\partial_i g),\partial_i g)\dd\\
 ={}&-\int_{S^d}\sum_iT_g(\Delta g,\partial_i g,\partial_i g)\dd
     -2\int_{S^d}\sum_i h_g(\partial_i\Delta g,\partial_i g)\dd.
\end{aligned}
\]

To integrate the last term by parts, consider the real one-form
\[
 \zeta(X):=h_g(\Delta g,\mathrm dg(X)),\qquad
 \zeta_i=\zeta(e_i)=h_g(\Delta g,\partial_i g).
\]
Here $\nabla_i(\Delta g)=\partial_i\Delta g$, since $\Delta g$ is a
matrix-valued function, while $(\nabla_i\mathrm dg)(e_j)=(\Hess g)_{ij}$. The covariant product rule
therefore gives
\[
 \operatorname{div}(\zeta^\sharp)=\sum_i(\nabla_i\zeta)_i
 =\sum_iT_g(\partial_i g,\Delta g,\partial_i g)
  +\sum_i h_g(\partial_i\Delta g,\partial_i g)+\sum_i h_g(\Delta g,(\Hess g)_{ii}).
\]
Applying divergence theorem~\eqref{eq:stokes} to the globally defined one-form $\zeta$,
and using $\sum_i (\Hess g)_{ii}=-\Delta g$ and the symmetry of $T_g$, yields
\[
 \int_{S^d}\sum_i h_g(\partial_i\Delta g,\partial_i g)\dd
 =\int_{S^d}\|\Delta g\|_g^2\dd
  -\int_{S^d}\sum_iT_g(\Delta g,\partial_i g,\partial_i g)\dd.
\]
Substituting this identity into the time derivative and
using~\eqref{eq:Tdiv}, we obtain
\begin{equation}\label{eq:J}
\begin{aligned}
 -\left. \dfrac12\dfrac{\mathrm d}{\mathrm dt}\I(g_t)\right|_{t=0^+}
 &=\int_{S^d}\|\Delta g\|_g^2\dd
   -\dfrac12\int_{S^d}\sum_iT_g(\Delta g,\partial_i g,\partial_i g)\dd\\
 &=\int_{S^d}\|\Delta g\|_g^2\dd-\mathcal C-\dfrac12\mathcal W.
\end{aligned}
\end{equation}

To retain the Hessian trace, recall that
\[
 (\Hess g)_{ij}=a^0_{ij}-\frac{\delta_{ij}}d\Delta g,\qquad
 \sum_i a^0_{ii}=0.
\]
The trace and trace-free parts are orthogonal under the inner product
induced by the spherical metric and $h_g$:
\[
 \sum_{i,j}h_g\left(a^0_{ij},\frac{\delta_{ij}}d\Delta g\right)
 =\frac1d h_g\left(\sum_i a^0_{ii},\Delta g\right)=0.
\]
Expanding the squared norm and integrating therefore gives
\begin{equation}\label{eq:Hessian-zero}
 \int_{S^d}\|\Hess g\|_g^2\dd
 =\int_{S^d}\sum_{i,j}\|a^0_{ij}\|_g^2\dd
  +\frac1d\int_{S^d}\|\Delta g\|_g^2\dd.
\end{equation}

For the mixed term, $\sum_i a^0_{ii}=0$ allows us to replace
$c^0_{ij}$ by $c_{ij}$. Using
$c_{ij}=c_g(\partial_i g,\partial_j g)$ and~\eqref{eq:c}, we obtain
\[
\begin{aligned}
 \sum_{i,j}h_g(a^0_{ij},c^0_{ij})
 &=\sum_{i,j}h_g(a^0_{ij},c_{ij})
  =\sum_{i,j}h_g((\Hess g)_{ij},c_{ij})+\frac1d\sum_i h_g(\Delta g,c_{ii})\\
 &=\sum_{i,j}T_g((\Hess g)_{ij},\partial_i g,\partial_j g)
   +\frac1d\sum_iT_g(\Delta g,\partial_i g,\partial_i g).
\end{aligned}
\]
Integrating and applying~\eqref{eq:Tdiv} yields
\begin{equation}\label{eq:interterm}
 \int_{S^d}\sum_{i,j}h_g(a^0_{ij},c^0_{ij})\dd
 =\mathcal C+\frac1d(2\mathcal C+\mathcal W)
 =\frac{d+2}{d}\mathcal C+\frac1d\mathcal W.
\end{equation}
Substituting~\eqref{eq:Hessian-zero} into the weighted Bochner
identity~\eqref{eq:weightedbochner}, we find
\[
 \frac{d-1}{d}\int_{S^d}\|\Delta g\|_g^2\dd
 =(d-1)\I(g)+\int_{S^d}\sum_{i,j}\|a^0_{ij}\|_g^2\dd
   +3\mathcal C+\mathcal W.
\]
Combining this with~\eqref{eq:J} gives
\[
 -\frac12\frac{\mathrm d}{\mathrm dt}\I(g_t)-d\I(g)
 =\frac d{d-1}\int_{S^d}\sum_{i,j}\|a^0_{ij}\|_g^2\dd
   +\frac{2d+1}{d-1}\mathcal C
   +\frac{d+1}{2(d-1)}\mathcal W.
\]
To eliminate the mixed term, expand the square and
use~\eqref{eq:interterm}:
\[
\begin{aligned}
 \int_{S^d}\sum_{i,j}\|a^0_{ij}+\beta_dc^0_{ij}\|_g^2\dd
 &=\int_{S^d}\sum_{i,j}
   \left(\|a^0_{ij}\|_g^2+2\beta_d h_g(a^0_{ij},c^0_{ij})
         +\beta_d^2\|c^0_{ij}\|_g^2\right)\dd\\
 &=\int_{S^d}\sum_{i,j}
   \left(\|a^0_{ij}\|_g^2+\beta_d^2\|c^0_{ij}\|_g^2\right)\dd
   +2\beta_d\left(\frac{d+2}{d}\mathcal C+\frac1d\mathcal W\right).
\end{aligned}
\]
Set $\displaystyle{\beta_d=\frac{2d+1}{2(d+2)}}$, the choice of $\beta_d$ cancels the $\mathcal C$ term after
completing the square, since
\[
 \frac{2d+1}{d-1}-2\beta_d\frac{d+2}{d-1}=0.
\]
The remaining coefficient of $\mathcal W$ is
\[
 \frac{d+1}{2(d-1)}-\frac{2\beta_d}{d-1}=\frac d{2(d+2)}.
\]
We thus obtain the exact identity
\begin{equation}\label{eq:exactcompletion}
 -\dfrac12\dfrac{\mathrm d}{\mathrm dt}\I(g_t)-d\I(g)
 =\dfrac d{d-1}\int_{S^d}\sum_{i,j}\|a^0_{ij}+\beta_dc^0_{ij}\|_g^2\dd
  +\dfrac d{2(d+2)}\mathcal W
  -\dfrac{d\beta_d^2}{d-1}\int_{S^d}\sum_{i,j}\|c^0_{ij}\|_g^2\dd.
\end{equation}
Lemma~\ref{lem:tensor}, with $k=d$ and $X_i=\partial_i g$, gives
\[
 \int_{S^d}\sum_{i,j}\|c^0_{ij}\|_g^2\dd
 \le\dfrac{d-1}{2d}\mathcal W.
\]
Hence the negative quadratic term in~\eqref{eq:exactcompletion}
is bounded below by
\[
 -\frac{d\beta_d^2}{d-1}\int_{S^d}\sum_{i,j}\|c^0_{ij}\|_g^2\dd
 \ge-\frac{\beta_d^2}{2}\mathcal W.
\]
This leaves a positive coefficient in front of $\mathcal W$, which
we denote by
\[
 \gamma_d:=\frac d{2(d+2)}-\frac{\beta_d^2}{2}
          =\frac{4d-1}{8(d+2)^2}>0.
\]
Moreover, $\mathcal W\ge0$ by~\eqref{eq:CW} and
Lemma~\ref{lem:differential}. Substituting the preceding bound into
\eqref{eq:exactcompletion} therefore yields
\[
 -\frac12\frac{\mathrm d}{\mathrm dt}\I(g_t)-d\I(g)
 \ge\frac d{d-1}\int_{S^d}\sum_{i,j}
       \|a^0_{ij}+\beta_dc^0_{ij}\|_g^2\dd+\gamma_d\mathcal W\ge0,
\]
which proves~\eqref{eq:bochnerlower}.
\end{proof}

\section{Applications to matrix concentration}
\label{sec:applications}

Theorem~\ref{thm:main} yields a matrix Laplace-transform bound through
a Herbst argument based on an exponential field with matrix mean.
Throughout this section,
expectations are taken with respect to normalized spherical measure on $S^d$.
For a matrix $A$, let $\|A\|_\infty$ and $\|A\|_1$ denote its operator
norm and trace norm, respectively.
For a smooth Hermitian field $F:S^d\to (M_n)_{\sa}$, set
\begin{equation*} \Gamma(F)=\sum_{j=1}^d(\partial_j F)^2,
\end{equation*}
where $\partial_j F=e_j(F)$ and $(e_j)$ is any local real orthonormal tangent frame. The resulting
positive semidefinite field $\Gamma(F)$ is independent of the choice of frame.
The matrix Laplace-transform approach of
Ahlswede--Winter~\cite[Appendix]{AW02} converts trace exponential
bounds into concentration estimates. Matrix entropy and Herbst
arguments are developed in~\cite[Section~7.2]{CT14}, and semigroup
methods in~\cite{HT21}. 

\begin{theorem}[Matrix concentration]\label{thm:app-concentration}
Let $d\ge2$, let $F:S^d\to (M_n)_{\sa}$ be smooth. Suppose that $\Gamma(F)(x)\preceq V$ for all $x$, where
$V\succeq 0$ is a constant matrix. Then, for every
$B=B^*\in M_n$ and every $\theta\in\mathbb R$,
\begin{equation}\label{eq:app-laplace}
 \E\tr e^{B+\theta (F-\E F)}
 \le \tr\exp\!\left(B+\frac{\theta^2}{2d}V\right).
\end{equation}
If $\|V\|_\infty>0$, then, for $t>0$,
\begin{align}
 \mu\{\lambda_{\max}(F-\E F)\ge t\}
 &\le n e^{-dt^2/(2\|V\|_\infty)},\label{eq:app-one-tail}\\
 \mu\{\|F-\E F\|_\infty\ge t\}
 &\le 2n e^{-dt^2/(2\|V\|_\infty)},\label{eq:app-norm-tail}\\
 \E\|F-\E F\|_\infty
 &\le\sqrt{\frac{2\|V\|_\infty\log(2n)}d}.\label{eq:app-expected-norm}
\end{align}
The coefficient $1/(2d)$ in
\eqref{eq:app-laplace} is optimal for all
matrix sizes.
\end{theorem}

\begin{proof}
Fix a faithful density matrix $\sigma$, so $\tr\sigma=1$.
We first bound $\theta\E\tr(G(F-\E F))-D(G\|\sigma)$ for smooth strictly positive
fields $G$ with $\E G=\sigma$. We will then apply this estimate to
a normalized exponential field to obtain the shifted Laplace bound.

For each $\theta\in\mathbb R$, we choose a constant Hermitian shift
$A_\theta$ so that $e^{A_\theta+\theta (F-\E F)}$ has mean $\sigma$.
To construct this shift, minimize
\[
 \mathcal F_\theta(A)
 =\E\tr e^{A+\theta (F-\E F)}-\tr(\sigma A)
 \qquad(A=A^*).
\]
Strict convexity follows from that of the trace exponential.
To prove coercivity, put
$u=\lambda_{\max}(A)$, $\ell=\lambda_{\min}(A)$,
$\eta=\lambda_{\min}(\sigma)>0$, and
$R=\sup_x\|F(x)-\E F\|_\infty$. Eigenvalue monotonicity and
$\tr(\sigma(uI-A))\ge\eta(u-\ell)$ give
\[
 \mathcal F_\theta(A)
 \ge e^{u-|\theta|R}-u+\eta(u-\ell).
\]
The term $e^{u-|\theta|R}-u$ tends to infinity as $|u|\to\infty$,
while $\eta(u-\ell)$ controls the spectral width. Thus
$\mathcal F_\theta$ is coercive and has a unique minimizer $A_\theta$.
For every Hermitian direction $K$, the first variation is
\[
 D_A\mathcal F_\theta(A)[K]
 =\tr\bigl((\E e^{A+\theta (F-\E F)}-\sigma)K\bigr).
\]
The first-order condition therefore gives the required fixed mean:
\begin{equation*} G_\theta=e^{A_\theta+\theta (F-\E F)},\qquad \E G_\theta=\sigma.
\end{equation*}
The Hessian is positive definite on the real space of Hermitian
matrices. The implicit function theorem therefore makes $A_\theta$
smooth in $\theta$. At $\theta=0$, $A_0=\log\sigma$.

To express the entropy of $G_\theta$, set
\mbox{$\psi_\sigma(\theta)
=\tr(\sigma\log\sigma)-\tr(\sigma A_\theta)$}.
Differentiating $\E\tr G_\theta=1$ and using $\E G_\theta=\sigma$ give
\[
 0=\frac{\mathrm{d}}{\mathrm{d}\theta}\E\tr G_\theta
  =\tr(\sigma A_\theta')+\E\tr(G_\theta (F-\E F)).
\]
Hence $\psi_\sigma'(\theta)=\E\tr(G_\theta (F-\E F))$. Since
$\log G_\theta=A_\theta+\theta (F-\E F)$, expanding the entropy yields
\[
 D(G_\theta\|\sigma)
 =\tr(\sigma A_\theta)+\theta\E\tr(G_\theta (F-\E F))
   -\tr(\sigma\log\sigma)
 =\theta\psi_\sigma'(\theta)-\psi_\sigma(\theta).
\]
At $\theta=0$, the identities $A_0=\log\sigma$ and $\E (F-\E F)=0$ give
$\psi_\sigma(0)=\psi_\sigma'(0)=0$.

We next bound the Fisher information. For every $G>0$ and Hermitian $K$,
\begin{equation*} \int_0^1\tr(G^sKG^{1-s}K)\,\mathrm{d}s\le\tr(GK^2).
\end{equation*}
In an eigenbasis of $G$, the left side is a sum weighted by scalar
logarithmic means. Bounding each by the arithmetic mean proves the
inequality. Since $A_\theta$ and $\E F$ are spatially constant,
$\partial_j\log G_\theta=\theta\partial_j F$. The gradient expression for
Fisher information and the derivative formula for the exponential now give
\begin{align*}
 \I(G_\theta)
 &=\theta^2\E\sum_j\int_0^1
 \tr\bigl(G_\theta^s(\partial_j F)G_\theta^{1-s}(\partial_j F)\bigr)\,\mathrm{d}s\\
 &\le\theta^2\E\tr(G_\theta\Gamma(F))
 \le\theta^2\tr(\sigma V).
\end{align*}
The last inequality uses $\Gamma(F)\preceq V$ and $\E G_\theta=\sigma$.
Theorem~\ref{thm:main} therefore gives
\[
 \theta\psi_\sigma'(\theta)-\psi_\sigma(\theta)
 \le\frac{\theta^2}{2d}\tr(\sigma V).
\]
For $\theta>0$, this is equivalent to
\[
 \left(\frac{\psi_\sigma(\theta)}{\theta}\right)'
 \le\frac1{2d}\tr(\sigma V).
\]
The initial conditions imply
$\displaystyle{\frac{\psi_\sigma(\theta)}{\theta}\to0}$
as $\theta\downarrow0$.
Integrating from zero to $\theta$ gives
\[
 \psi_\sigma(\theta)\le\frac{\theta^2}{2d}\tr(\sigma V).
\]
Applying the same argument to $-F$ proves this bound for negative
$\theta$ as well.

To pass from $G_\theta$ to an arbitrary smooth strictly positive field $G$
with $\E G=\sigma$, we use Klein's inequality:
\begin{equation}\label{eq:app-variational}
 \theta\E\tr(G(F-\E F))-D(G\|\sigma)\le\psi_\sigma(\theta).
\end{equation}
Indeed, the difference between the right and left sides is
\[
 \E\tr\bigl(G(\log G-\log G_\theta)\bigr)\ge0,
\]
because $\E\tr G=\E\tr G_\theta=1$.
Combining~\eqref{eq:app-variational} with the quadratic bound for
$\psi_\sigma$ yields the prescribed-mean estimate
\begin{equation}\label{eq:app-constrained-herbst}
 \theta\E\tr(G(F-\E F))-D(G\|\sigma)
 \le\frac{\theta^2}{2d}\tr(\sigma V),
 \qquad \E G=\sigma,\quad\tr\sigma=1.
\end{equation}
Since~\eqref{eq:app-constrained-herbst} holds for every faithful density
matrix $\sigma$, we apply it with
\[
 Z=\E\tr e^{B+\theta (F-\E F)},\qquad
 G=Z^{-1}e^{B+\theta (F-\E F)},\qquad \sigma=\E G.
\]
Then $\sigma$ is a faithful density matrix and
$\log G=B+\theta (F-\E F)-(\log Z)I_n$. Expanding the entropy and
applying~\eqref{eq:app-constrained-herbst}, we obtain
\begin{align*}
 \log Z
 &=\theta\E\tr(G(F-\E F))-D(G\|\sigma)
   +\tr(\sigma B)-\tr(\sigma\log\sigma)\\
 &\le\tr\!\left(\sigma\left(B+\frac{\theta^2}{2d}V\right)\right)
       -\tr(\sigma\log\sigma)\\
 &\le\log\tr\exp\!\left(B+\frac{\theta^2}{2d}V\right).
\end{align*}
The last inequality follows by comparing $\sigma$ with the normalized
exponential of $\displaystyle{B+\frac{\theta^2}{2d}V}$
in the Gibbs variational principle.
This proves~\eqref{eq:app-laplace}.

For the tail bounds, take $B=0$ to obtain
$\E\tr e^{\theta (F-\E F)}\le n e^{\theta^2\|V\|_\infty/(2d)}$.
For $\theta>0$, the exponential Markov inequality~\cite[Proposition~3.1]{Tropp12} gives
\[
 \mu\{\lambda_{\max}(F-\E F)\ge t\}
 \le e^{-\theta t}\E\tr e^{\theta (F-\E F)}
 \le n\exp\!\left(-\theta t+\frac{\theta^2\|V\|_\infty}{2d}\right).
\]
The choice $\displaystyle{\theta=\frac{dt}{\|V\|_\infty}}$
gives~\eqref{eq:app-one-tail}.
Applying the same estimate to $-F$ and taking
a union bound give~\eqref{eq:app-norm-tail}. For the expected norm, use
\[
 e^{\theta\|F-\E F\|_\infty}\le\tr e^{\theta (F-\E F)}+\tr e^{-\theta (F-\E F)}
 \quad(\theta>0).
\]
Jensen's inequality then yields
$\displaystyle{\E\|F-\E F\|_\infty
\le\frac{\log(2n)}{\theta}+\frac{\theta \|V\|_\infty}{2d}}$.
Optimizing proves~\eqref{eq:app-expected-norm}. If $V=0$, the gradient
vanishes and $F$ is constant.

To prove optimality, choose finite-dimensional Hermitian Clifford matrices
$C_1,\ldots,C_{d+1}$ with
$C_iC_j+C_jC_i=2\delta_{ij}I_n$, and put
\[
 F(x)=\frac1{\sqrt d}\sum_{j=1}^{d+1}x_jC_j.
\]
These Clifford matrices are traceless, $\E F=0$, and
$\displaystyle{F(x)^2=\frac{I_n}{d}}$.
Using $\langle\nabla x_i,\nabla x_j\rangle=\delta_{ij}-x_ix_j$,
\[
 \Gamma(F)=\frac1d\sum_j C_j^2-F^2=I_n,
 \qquad
 \frac1n\E\tr e^{\theta F}
 =\cosh\!\left(\frac{\theta}{\sqrt d}\right).
\]
Since
$\displaystyle{\log\cosh\!\left(\frac{\theta}{\sqrt d}\right)
 =\frac{\theta^2}{2d}+O(\theta^4)}$,
no smaller universal quadratic coefficient is possible.
\end{proof}

The Ricci curvature estimate of
Huang--Tropp~\cite[Theorem~3.2 and Section~4.3.1]{HT21} gives the
denominator $d-1$ on $S^d$ under the same scalar gradient bound.
Theorem~\ref{thm:app-concentration} gives $d$, with a matrix bound $V$
and an arbitrary Hermitian shift $B$. Under $\Gamma(F)\preceq I_n$,
the quadratic trace-Laplace coefficient is optimal; the sharpness
assertion does not concern the tail prefactors.

For example, let $F(x)=\sum_{j=1}^{d+1}x_jA_j$ for Hermitian matrices
$A_1,\ldots,A_{d+1}$, and put $S=\sum_jA_j^2$. Then
\[
 \E F=0,\qquad \Gamma(F)=S-F^2\preceq S,\qquad
 \E\tr e^{B+\theta F}\le\tr e^{B+\theta^2S/(2d)}.
\]
Thus the bound applies to matrix series with dependent spherical
coefficients.

\begin{samepage}
\begin{example}[Rank-one sensing and spherical designs]
Let $n\ge2$ and let $a_1,\ldots,a_m\in S^{n-1}$ form a real spherical $2$-design,
so finite averages reproduce sphere averages for polynomials of
degree at most two~\cite[Section~5]{DGS77}. Put
\[
 A_j=a_ja_j^\top,\qquad A_j^2=A_j,\qquad
 \sum_{j=1}^mA_j=\frac mn I_n.
\]
For example, the union of two orthonormal bases in general position
and their antipodes gives $m=4n$ directions.
For a real symmetric matrix $X$, define
$\mathcal A(X)_j=\tr(A_jX)$. For $0\ne u\in\mathbb R^n$, the rank-one signal
$X_\star=uu^\top$ gives the quadratic measurements
$\mathcal A(X_\star)_j=|\langle a_j,u\rangle|^2$, as in the
PhaseLift formulation of phase retrieval~\cite{CSV13}.

Consider noisy observations
$y=\mathcal A(X_\star)+\sigma\xi$, where $\sigma>0$ and
$\xi$ is uniform on $S^{m-1}$, with $m\ge3$.
Thus the noise has fixed Euclidean norm $\sigma$.
The adjoint noise is the real symmetric matrix
\[
 Z=\mathcal A^*(y-\mathcal A(X_\star))
   =\sigma\sum_{j=1}^m\xi_jA_j.
\]
The concentration variable is the noise direction $\xi\in S^{m-1}$.
The linear-series calculation above gives
\[
 \E Z=0,\qquad
 \Gamma(Z)=\sigma^2\sum_jA_j^2-Z^2
          =\frac{\sigma^2m}{n}I_n-Z^2
          \preceq\frac{\sigma^2m}{n}I_n.
\]
Applying Theorem~\ref{thm:app-concentration} with $d=m-1$ yields,
for $0<\delta<1$,
\[
 \mathbb P\!\left\{\|Z\|_\infty\ge
 \sigma\sqrt{\frac{2m}{n(m-1)}\log\frac{2n}{\delta}}\right\}
 \le\delta.
\]
This gives an explicit noise threshold for trace-norm regularization
of the lifted signal. Since $\|X\|_1=\tr X$ on the positive cone,
consider, for example,
\[
 \widehat X_\lambda\in\operatorname*{arg\,min}_{X=X^\top\succeq0}
 \left(\frac12\|\mathcal A(X)-y\|_2^2+\lambda\tr X\right).
\]
Taking $\lambda$ to be twice the noise threshold above ensures
$\lambda\ge2\|\mathcal A^*(y-\mathcal A(X_\star))\|_\infty$
with probability at least $1-\delta$. On this event, the stochastic term
$\tr(Z(\widehat X_\lambda-X_\star))$ in the least-squares basic inequality
is bounded in absolute value by
$(\lambda/2)\|\widehat X_\lambda-X_\star\|_1$.
\end{example}
\end{samepage}

\bibliographystyle{amsplain}
{\raggedright
\bibliography{sharp_complete_entropy_draft}
}

\end{document}